\documentstyle[amsfonts]{article}
\begin{document}
\newtheorem{theor}{Theorem}[section] 
\newtheorem{prop}[theor]{Proposition} 
\newtheorem{cor}[theor]{Corollary}
\newtheorem{lemma}[theor]{Lemma}
\newtheorem{sublem}[theor]{Sublemma}
\newtheorem{defin}[theor]{Definition}
\newtheorem{conj}[theor]{Conjecture}

\hfuzz2cm

\gdef\Aut{{\rm Aut}}
\gdef\beginProof{\par{\bf Proof: }}
\gdef\endProof{${\bf Q.E.D.}$\par}  
\gdef\ar#1{\widehat{#1}}
\gdef\pr{^{\prime}}   
\gdef\prpr{^{\prime\prime}}
\gdef\mtr#1{\overline{#1}}
\gdef\ra{\rightarrow}
\gdef\Bbb{\bf }
\gdef\a{\alpha}
\gdef\ca{ch_{g}(\alpha)}
\gdef\rl{{\Lambda}}
\def\covol{{\rm covol}}
\def\lcovol{{\rm lcovol}}
\gdef\CT{CT}
\gdef\refeq#1{(\ref{#1})}
\gdef\mn{{\mu_{n}}}
\gdef\zn{{\Bbb Z}/(n)}
\gdef\umn{^{\mn}}
\gdef\lmn{_{\mn}}
\gdef\blb{{\big(}}
\gdef\brb{{\big)}}
\gdef\ttmk#1{\widetilde{\lambda}(#1)}
\gdef\chge#1{ch_{g}^{-1}(\lambda_{-1}({#1}^{\vee}))
ch_{g}(\lambda_{-1}({#1}^{\vee}_{g}))}
\gdef\mlcovol{{\rm mlc}}
\gdef\Hom{{\rm Hom}}
\def\Trs{{\rm Tr_s\,}}
\def\Tr{{\rm Tr\,}}
\def\End{{\rm End}}
\def\eq{equivariant }
\def\Td{{\rm Td}}
\def\ch{{\rm ch}}
\def\torus{{\cal T}}
\def\Proj{{\rm Proj}}
\def\bmn{{R_{n}}}
\def\uexp#1{{{\rm e}^{#1}}}
\def\Spec{{\rm Spec}\,}
\def\Qb{\mtr{\Bbb Q}}
\def\Cn{{\bf C}_n}
\def\Zn{{{\Bbb Z}/n}}
\def\deg{{\rm deg}}
\def\mod{{\rm mod}}
\def\ac1{\ar{\rm c}_{1}}
\def\boxtimes{{\otimes_{\rm Ext}}}
\def\Qmn{{{\Bbb Q}(\mu_{n})}}
\def\NIm#1{{\rm Im}(#1)}
\def\NRe#1{{\rm Re}(#1)}
\def\rk{{\rm rk}\,}
\def\rrot{{{\rm R}^{\rm rot}}}
\newcommand{\kkk}[1]{{\Large\bf#1}}
\newcommand{\nlabel}[1]{\label{#1}}

\author{ Kai K\"ohler
\footnote{Mathematisches Institut, Wegelerstr. 10, 
D-53115 Bonn, Germany, E-mail: koehler@math.uni-bonn.de}\\
         Damian Roessler
\footnote{Centre de Math\'ematiques de Jussieu,  
     Universit\'e Paris 7 Denis Diderot,  
     Case Postale 7012,  
     2, place Jussieu,  
     F-75251 Paris Cedex 05, France,   
     E-mail : roessler@math.jussieu.fr}\\}
\title{A fixed point formula of Lefschetz type in Arakelov geometry IV: 
the modular height of C.M. abelian varieties}
\maketitle
\begin{abstract}
We give a new proof of a slightly weaker form of a theorem 
of P. Colmez (\cite[Par. 2]{C2}). This theorem 
(Corollary \ref{MainCor}) gives a formula for 
the Faltings height of abelian varieties 
with complex multiplication by a C.M. field whose 
Galois group over $\bf Q$ is abelian; it reduces 
to the formula of Chowla and Selberg in the case of elliptic 
curves. We show that the 
formula can be deduced from the arithmetic fixed point formula 
 proved in \cite{KR2}. Our proof is intrinsic in the sense that 
it does not rely 
on the computation of the periods of any particular abelian 
variety.
\end{abstract}
\begin{center}
\end{center}
\begin{center}
\end{center}
\thispagestyle{empty}
\newpage
\setcounter{page}{1}
\tableofcontents
\newpage
\parindent=0pt
\parskip=5pt

\section{Introduction}

Let $A$ be an abelian variety of dimension $d$ defined over $\mtr{\bf Q}$. 
Let $K\subseteq {\bf C}$ be a number field such that $A$ is defined over $K$ and such 
that the N\'eron model of 
$A$ over 
${\cal O}_{K}$ has semi-stable reduction at all the places of $K$.
 Let $\Omega$ be the ${\cal O}_{K}$-module 
of global sections of the sheaf of differentials of ${\cal A}$ over 
${\cal O}_{K}$ and let $\alpha$ be a section of $\Omega^{d}$. 
We write ${A}({\bf C})_{\sigma}$ for the manifold of complex points 
of the variety ${A}\times_{\sigma({K})}{\bf C}$, where 
$\sigma\in\Hom(K,{\bf C})$ is an embedding of ${K}$ in $\bf C$. 
The
modular (or Faltings) height of $A$ is the quantity
$$
h_{\rm Fal}(A):={1\over[K:{\bf Q}]}\log(\#\Omega^{d}/\alpha.\Omega^{d})-
{1\over 2[K:{\bf Q}]}\sum_{\sigma:K\ra{\bf C}}\log|{1\over(2\pi)^{d}}
\int_{{\cal A}({\bf C})_{\sigma}}
\alpha\wedge\mtr{\alpha}|.
$$
It does not depend on the choice of $K$ or $\alpha$. 
The modular height defines a height on some moduli spaces of 
abelian varieties and plays a key role in Falting's proof of 
the Mordell conjecture. The object of this article 
is to use higher dimensional Arakelov theory to prove a formula 
for the modular height of $A$, valid if 
$A$ has complex multiplication by the ring of integers 
of an abelian extension of $\bf Q$. A full self-contained 
statement of the formula can in be found in 
corollary \ref{MainCor}. This formula was first proved 
by a completly different method by P. Colmez and the 
remaining of the introduction is devoted to 
an exposition of his and our approach, followed 
by a plan of the paper. \\
If $A$ is an abelian variety with complex multiplication, 
the modular height is 
related to the periods of $A$. So suppose that there exists 
a C.M. field $E$, of degree $2d$ over $\bf Q$ and an 
embedding of rings ${\cal O}_{E}\ra {\rm End}(A)$ into the 
endomorphism ring of $A$. This is equivalent to 
saying that $A$ has complex multiplication  by ${\cal O}_{E}$ 
(cf. \cite{Sh}).  
We can suppose without loss of generality that the action 
of ${\cal O}_E$ is defined over $K$ and that $K$ contains all 
the conjugates of $E$ in $\bf C$. Let now $\tau\in\Hom(E,{\bf C})$ be 
an embedding of $E$ in $\bf C$ and let $\omega_{\tau}$ be 
an element of the first algebraic de Rahm cohomology group 
$H^1_{\rm DR}(A)$ of $A$ over $K$ (this is a $K$-vector space of rank $2d$), such 
that $a(\omega_{\tau})=\tau(a).\omega_\tau$ for all 
$a\in{\cal O}_{E}$. If $\sigma\in\Hom(K,{\bf C})$, 
let $\omega_{\tau}^{\sigma}$ be the element of the complex 
cohomology group $H^{1}(A({\bf C})_{\sigma},{\bf C})$ obtained by 
base change. Let $u_{\sigma}$ be an element 
of the rational homology group $H_{1}(A({\bf C})_\sigma, {\bf Q})$. 
The {\it period} $P(A,\tau,\sigma)\in{\bf C}/K^*$ 
associated to $\tau$ and $\sigma\in 
\Hom(K,{\bf C})$ is the complex 
number $<\omega_{\tau}^\sigma,u_{\sigma}>_\infty$, where 
$<\cdot,\cdot>_\infty$ is the natural pairing between cohomology 
and homology. Up to multiplication by an element of $K^{*}$, it 
is only dependent on $\tau$ and $\sigma$. 
Let $\Phi$ be the subset of $\Hom(E,{\bf C})$ such that 
the subspace $\{t\in TA_{0}:a(t)=\tau(a).t,\ \forall a\in
{\cal O}_E\}$ contains 
a non-zero element (this is the {\it type} of the C.M. abelian 
variety $A$). 
The 
following lemma is a (very) weak form of a theorem of 
P. Colmez: 
\begin{lemma}
The equality
$$
(2\pi)^{-d/2}.e^{-[K:{\bf Q}]h_{\rm Fal}(A)}=\prod_{\tau\in\Phi}\prod_{\sigma\in
\Hom(K,{\bf C})}P(A,\tau,\sigma)
$$
holds up to multiplication by an element of $\mtr{\bf Q}$. 
\end{lemma}
Furthermore, using a refinement of the above 
lemma, the theory of $p$-adic periods and 
explicit computation of periods of Jacobians of Fermat 
curves, he gives an explicit formula for $h_{\rm Fal}(A)$ 
(see also \cite{And},  
\cite{Yo} and \cite{Gr} [for elliptic curves] 
for mod. $\mtr{\bf Q}$ versions 
of the latter formula). 
To describe it, suppose furthermore that 
$E$ is Galois over $\bf Q$ and let $G:={\rm Gal}(E|{\bf Q})$. 
Identify 
$\Phi$ with its characteristic function $G\ra\{0,1\}$ and  
define $\Phi^{\vee}$ by the formula $\Phi^{\vee}(\tau):=
\Phi(\tau^{-1})$. 
\begin{theor}[Colmez]
If $G$ is abelian, there exists $q\in{\bf Q}$, such that the identity
$$
{1\over d}h_{\rm Fal}(A)=-\sum_{\chi\ {\rm odd}}<\Phi*\Phi^{\vee},\chi>
[2{L\pr(\chi_{\rm prim},0)\over L(\chi_{\rm prim},0)}+
\log(f_{\chi})]+q\log(2)
$$
holds. If the conductor of $E$ over $\bf Q$ 
divides $8n$, where $n$ is an odd natural number, then 
the identity holds with $q=0$. 
\nlabel{MColTh}
\end{theor}
(it is conjectured that $q$ always vanishes) 
Here $<\cdot,\cdot>$ refers to the scalar product 
of complex valued functions on $G$ and $*$ to the 
convolution product. The sum $\sum_{\chi\ {\rm odd}}$ is 
on all the odd characters of $G$ (recall that 
$\chi$ is odd iff $\chi(h\circ c\circ h^{-1}\circ\tau)=-\chi(\tau)$ for all
$\tau,h\in G$, where $c\in G$ is complex conjugation). The 
notation $f_{\chi}$ refers to the conductor of $\chi$.
Colmez conjectures that the theorem \ref{MColTh} holds even 
without the condition that $G$ is abelian. This formula can be viewed 
as a generalisation of the formula of Chowla and Selberg 
(see \cite{CS}), 
to which it reduces when applied to a C.M. elliptic curve.\\
It is the aim of this paper to provide a proof of 
\ref{MColTh} using higher dimensional Arakelov theory. More 
precisely, we shall show that a slightly weaker form 
of \ref{MColTh} can be derived from the fixed point 
formula in Arakelov theory proved in \cite{KR2} (announced 
in \cite{KR1}), when applied 
to abelian varieties with complex multiplication
by a field generated over $\bf Q$ by a root of $1$.
 This proof has the advantage 
of being intrinsic, i.e. the right side 
of \ref{MColTh} is obtained directly from analytic invariants 
(the equivariant analytic torsion and the equivariant $R$-genus) 
of the abelian variety. It does not involve the 
computation of the periods of a particular C.M. abelian 
variety (e.g. Jacobians of Fermat curves). We shall prove:
\begin{theor}
Let $f$ be the conductor of $E$ over $\bf Q$; 
if $G$ is abelian 
there exist numbers $a_{p}\in{\bf Q}(\mu_f)$, where $p|f$, such 
that the identity
$$
{1\over d}h_{\rm Fal}(A)=-\sum_{\chi\ {\rm odd}}<\Phi*\Phi^{\vee},\chi>
2{L\pr(\chi_{\rm prim},0)\over L(\chi_{\rm prim},0)}+
\sum_{p|f}a_{p}\log(p)
$$
holds. 
\nlabel{MainTh}
\end{theor}
Notice that the difference between the right sides of 
the equations in \ref{MColTh} and \ref{MainTh} is equal to 
$\sum_{p|f}b_{p}\log(p)$ for some $b_{p}\in{\bf Q}(\mu_f)$. 
To see this, let us write 
$$
\sum_{\chi\ {\rm odd}}<\Phi*\Phi^{\vee},\chi>\log(f_{\chi})=
\sum_{p|f}[\sum_{\chi\ {\rm odd}}<\Phi*\Phi^{\vee},\chi>n_{\chi,p}]\log(p)
$$
where $n_{\chi,p}$ is the multiplicity of the prime number 
$p$ in the number $f_{\chi}$. By construction, the number 
$\sum_{\chi\ {\rm odd}}<\Phi*\Phi^{\vee},\chi>n_{\chi,p}$ is invariant 
under the action of ${\rm Gal}(\mtr{\bf Q}|{\bf Q})$ and is thus 
rational.\\
The paper is organised as follows. In section 2, we 
give the rephrasing of the main result of \cite{KR2} which 
we shall need in the text. In section 3, we prove 
an equivariant version of Zarhin's trick (this is 
of independent interest) and we show 
that the fixed point scheme of the action of a finite 
group on an abelian scheme is well-behaved away 
from the fibers lying over primes dividing 
the order of the group. In section 
4 we compute a term occuring in the fixed point formula, namely the
equivariant holomorphic torsion of line bundles on complex tori; we
adapt a method by Berthomieu to do so. In section 5, we apply 
the fixed point formula to the following setting: 
an abelian scheme with the action of a certain group 
of roots of unity and an equivariant ample 
vector bundle with Euler characteristic $1$, which 
is provided by section 3; an expression for the Faltings 
height is then very quickly obtained (as a solution of 
the system of equations \refeq{MainEq}) and the rest of the section 
is concerned with equating this expression with the linear 
combination of logarithmic derivatives of $L$-functions 
appearing in \ref{MColTh}.\\
We shall in this paper freely use the definitions 
and terminology of section 4 of \cite{KR2} (mainly up to 
def. 4.1 included).\\ 
{\bf Acknowledgments.} Our thanks go mainly to 
P. Colmez for many interesting conversations and hints 
and for kindly providing a proof of corollary \ref{ColLem}. It is also 
a pleasure to thank V. Maillot and C. 
Soul\'e for stimulating discussions. Thanks as well 
to G. Faltings for pointing out a mistake 
in an earlier version of this text (a hypothesis 
was missing in lemma \ref{FaltLem}) and V. Maillot again for pointing 
out a redundance. We thank the SFB 256, 
"Nonlinear Partial Differential
Equations", at the University of Bonn for its support. The second author 
is grateful to the IHES (Bures-sur-Yvette) for its support in 
1998-99.

\section{An arithmetic fixed point formula}

In this section, we formulate the 
fixed point formula we shall apply to 
abelian varieties (more precisely, to their N\'eron models). 
It is an immediate consequence of the more 
general fixed point formula which is the main result 
\cite[Th. 4.4]{KR2} of \cite{KR2}. 
To formulate it, we shall 
need the 
notions of arithmetic Chow theory, arithmetic degree and 
arithmetic characteristic classes; for these see \cite[Par. 2.1]{BoGS}. 
To understand how it can be deduced from 
\cite[Th. 4.4]{KR2}, one 
needs to read the last section of \cite{KR2}. 
We shall nevertheless give a self-contained presentation (not proof) 
of the formula, when the scheme is smooth, and 
the associated fixed point scheme \'etale over the base scheme.\\
Let now $K$ be a number field and ${\cal O}_{K}$ its ring of 
integers. Let $U$ be an open subset of ${\rm Spec}{\bf Z}$ and 
let $U_{K}$ be the open subset of ${\rm Spec}{\cal O}_{K}$ lying 
above $U$. Let $p_{1},\dots,p_{l}$ be the prime numbers in the complement 
of $U$. Notice that $U_{K}$ is the spectrum of an arithmetic 
ring (see the beginning of \cite[Sec. 4]{KR2} for the definition). 
Denote by $\ar{\deg}:\ar{\rm CH}^{1}(U_{K})\ra\ar{\rm CH}^{1}(U)$ 
the push-forward map in arithmetic Chow theory (this map 
coincides with the usual arithmetic degree map if 
$U_{K}=\Spec{\cal O}_{K}$). By abuse of language, we shall 
write $\ar{\deg}(\mtr{V})$ for the arithmetic 
degree $\ar{\deg}(\ac1(\mtr{V}))$ of the first arithmetic 
Chern class of a hermitian bundle $\mtr{V}$ on $U_K$. Recall 
that the formula 
$$
\ar{\deg}(\ac1(\mtr{V}))=
\ar{\deg}(\ac1({\rm det}(\mtr{V})))=
$$
$$
\log(\#({\rm det}(V)/s.U_{K}))-
\sum_{\sigma\in\Hom(K,{\bf C})}\log|s\otimes_{\sigma(U_{K})}1|_{\sigma}
$$
holds. Here $s$ is any section of ${\rm det}(V)$ and 
$|\cdot|_{\sigma}$ is the norm arising from 
 the hermitian metric of  
${\rm det}(V)\otimes_{\sigma(U_{K})}{\bf C}$. 
It is a consequence of 
\cite[I, Prop. 2.2, p. 13]{GS3} that there is a canonical isomorphism 
$\ar{\rm CH}^{1}(U)\simeq {\Bbb R}/{\cal R}$, where 
$\cal R$ is the subgroup of $\Bbb R$ consisting of the 
expressions $\sum_{k=1}^{l}n_{k}\log p_{k}$, where 
$n_k\in{\Bbb Z}$. We shall thus often identify both.\\ 
Let now $n\in {\bf N}^{*}$ and let 
$\cal C$ be the subgroup of $\Bbb C$ generated by
 the expressions
$\sum_{k=1}^{l}q_{k}\alpha_{k}\log p_{k}$, where 
$q_{k}\in{\Bbb Q}$ is a rational number and 
$\alpha_{k}$ is an $n$-th 
root of unity. Let $f:Y\ra U_{K}$ be a
$\mn$-equivariant arithmetic variety 
over $U_{K}$, of relative dimension $d$. As usual, we fix 
a primitive root of unity $\zeta_{n}$. 
Recall also that $g$ is the automorphism of the complex manifold 
 $Y({\bf C})$ associated
to $\zeta_{n}$ via the action of $\mn({\bf C})$ on $Y({\bf C})$.
 Fix a $g$-invariant K\"ahler metric on $Y({\bf C})$, with 
associated K\"ahler form $\omega_{Y}$. 
\begin{theor}
Let $\mtr{E}$ be a $\mn$-equivariant hermitian vector bundle on $Y$. 
Suppose that $R^{k}f_{*}E=0$ for $k>0$. Then the equality 
\begin{eqnarray*}\lefteqn{
\sum_{k\in{\Bbb Z}/(n)}{\zeta_{n}^k}.\ar{\deg}((R^{0}f_{*}\mtr{E})_k)=
}\\&=&
{1\over 2}T_{g}(Y({\Bbb C}),\mtr{E})-
{1\over 2}\int_{Y_{\mu_{n}}({\bf C})}\Td_{g}(TY)\ch_g(E)R_{g}
(TY_{\bf C})+
\ar{{\rm deg}}\blb f_{*}\umn(\ar{\Td}_{\mn}(\mtr{Tf})
\ar{\ch}_{\mu_{n}}(\mtr{E}))\brb
\end{eqnarray*}
holds in ${\bf C}/{\cal C}$.
\nlabel{RBisConj}
\end{theor}
Recall that $R^{0}f_{*}\mtr{E}$ refers to the $U_{K}$-module 
of global  sections of $E$, endowed with the $L_{2}$-metric 
inherited from the hermitian metric on $E$ and the 
K\"ahler metric on $Y({\bf C})$. Recall that the 
$L_{2}$-metric is defined as follows:  
if $s,l$ are two holomorphic sections of $E_{\bf C}$, then 
$$
<s,l>_{L_{2}}:={1\over d!(2\pi)^{d}}\int_{Y({\bf C})}<s,l>_{E}
\omega_{Y}^{d}.
$$ 
\beginProof
The proof is similar to the proof of \cite[Th. 6.14]{KR2} and 
so we omit it.
\endProof
The above theorem is an immediate consequence of the arithmetic 
Riemann-Roch theorem of Bismut-Gillet-Soul\'e when the 
equivariant structure is trivial. 
Suppose now that $Y$ is smooth and $Y\lmn$ \'etale 
over $U_{K}$. Let 
$$
{L}^{\rm Im}(z,s):=\sum_{k\geq 0}{{\rm Im}(z^k)\over k^{s}}.
$$
where $z\in{\bf C}$, $|z|=1$ and $s\in{\bf C}$, $s>1$. 
As a function of $s$, it extends to a meromorphic function of 
the whole plane, which is holomorphic at $s=0$. 
We write $\rrot(\arg(z))$ for the derivative 
${\partial\over \partial s}{L}^{\rm Im}(z,0)$ of 
the function ${L}^{\rm Im}(z,s)$ at $0$. 
Write  
$\zeta$ for $\zeta_{n}$. We let $\Omega$ be the sheaf of differentials of 
$f$, which is locally free. Define
$$
D:=\prod_{k\in\Zn}(1-\zeta^k)^{\rk(\Omega_k)},\ 
T=\sum_{k\in\Zn}\zeta^{k}\rk(E_{k})
$$
which are locally constant complex-valued functions on $Y\lmn$.   
In the just mentionned setting, the formula in theorem 
\ref{RBisConj} becomes
\begin{eqnarray}\lefteqn{
\sum_{k\in{\Bbb Z}/(n)}{\zeta_{n}^k}.\ar{\deg}((R^{0}f_{*}\mtr{E})_k)=
}\nonumber\\&=&
{1\over 2}T_{g}(Y({\Bbb C}),\mtr{E})-
i\sum_{P\in Y_{\mu_{n}}({\bf C})}{{\rm Trace}(g|_{E_{P}})\over 
{\rm Det}({\rm Id}-g|_{\Omega_{P}})}.[\sum_{k\in\Zn}
\rk((TY({\bf C})_{P})_{k})
\rrot(\arg(\zeta^{k}))
]\nonumber\\
&+&
\ar{\deg}\blb f\umn_{*}({1\over D}\sum_{k\in\Zn}\zeta^{k}\ac1(\mtr{E}_{k})+
{T\over D}\sum_{k\in\Zn}{\zeta^{k}\over 1-\zeta^{k}}
\ac1(\mtr{\Omega}_{k}))\nlabel{RBCEq}\brb
\end{eqnarray}
in ${\bf C}/{\cal C}$. 
Finally, let us recall the following lemma.
\begin{lemma}
Let $n\geq 1$ be a a natural number and let $R$ be an entire ring 
of characteristic $0$ in which 
the polynomial $X^{n}-1$ splits. Let $\Cn$ be the constant group scheme over 
$\Bbb Z$ associated to $\Zn$, the cyclic group of order $n$. For each 
primitive $n$-th root of unity in $R$, there is 
an isomorphism of group schemes
$$
\mn\times_{\Bbb Z}\Spec R[{1\over n}]
\simeq \Cn\times_{\Bbb Z}\Spec R[{1\over n}]
$$
over $R[{1\over n}]$. 
\nlabel{CNMN}
\end{lemma}
\beginProof Both group 
schemes define isomorphic sheaves on the (small)
\'etale site of $R[{1\over n}]$ and they are both 
\'etale over $R[{1\over n}]$ (see for instance 
\cite[Par. 3.1, p. 100]{Tam}). Hence 
they are isomorphic.
\endProof
Let as usual $\Qmn$ refer to the number field generated 
by the complex $n$-th roots of unity. 
In view of the lemma, the constant group scheme $\Cn$ and 
the group scheme $\mn$ become isomorphic over $\Qmn[{1\over n}]$.
If we let $U$ be the set of prime divisors of $n$, then 
$\Qmn[{1\over n}]$ corresponds  to the open set $U_\Qmn$. Hence 
we see that the formula above can be applied to the action 
of an automorphism of finite order of a (regular, integral, 
projective) scheme over $\Qmn[{1\over n}]$


\section{Equivariant geometry on abelian schemes}

The following two lemmata give an equivariant version of Zarhin's trick. 
Let $\mtr{K}$ be an algebraically closed field of characteristic zero. 
If $X$ is a $\Zn$-equivariant abelian variety over $\mtr{K}$ and 
$L$ is a $\Zn$-equivariant line bundle on $X$, we shall 
say that the action of $\Zn$ on $L$ is {\it normalised}, if 
it induces the trivial action on the fiber $L|_{0}$ of $L$ 
at the origin. In this section, we shall write 
$a(\cdot)$ for the action of $1\in\Zn$. 
\begin{lemma}
Let $X,Y$ be $\Zn$-equivariant abelian varieties defined over $\mtr{K}$.
 Suppose 
that $\Zn$ acts by isogenies on $X$ and $Y$. Let $f:X\ra Y$ be 
an equivariant isogeny and let $N$ be the kernel
of $f$. Let $M$ be a line bundle on $X$. Suppose that $M$ is 
endowed with an $N$-equivariant structure and with a normalised
$\Zn$-equivariant structure. Let 
$\alpha\in \Hom(N_\Zn,\mtr{K}^{*})
\simeq H^{1}(N_\Zn,\mtr{K}^{*})$ be defined by the 
formula $\alpha(n)=a^{-1}\circ n\circ a\circ (-n)$. 
Then there exists a normalised $\Zn$-equivariant line bundle $L\pr$ on $Y$ 
and a $\Zn$-equivariant isomorphism $f^{*}L\pr\simeq L$, if and 
only if $\alpha=1$.
\nlabel{FaltLem}
\end{lemma}
Note that we have used the identification $\Aut(M)\simeq 
\mtr{K}^{*}$ in our definition of $\alpha$. Recall also 
that $N_\Zn$ refers to the part of $N$ fixed by every element 
of $\Zn$. 
\beginProof
If there exists a bundle $L\pr$ satisfying the 
hypothesies of the  lemma, then the equivariant structure 
of $f^{*}L\pr|_{N_\Zn}$ is by construction trivial
and thus $\alpha=1$. So suppose that $\alpha=1$. 
Note that $N$ 
is sent into itself by the elements of $\Zn$. 
Let $\rho$ be the character of 
$N$ defined by the formula 
$a^{-1}\circ a(n)\circ a\circ (-n)$ (this character extends 
$\alpha$). Since $\alpha=1$, 
$\rho$ induces a character $\rho\pr$ on the quotient group 
$N/N_\Zn$. This quotient is naturally identified 
with the image of the endomorphism $a-{\rm Id}$ of $N$. View $\rho\pr$ as 
a character on ${\rm Im}(a-{\rm Id})$ and choose any character $\rho\prpr$ 
extending $(\rho\pr)^{-1}$ to $N$. Such a character always exists because 
$\mtr{K}^{*}$ is an injective abelian group. 
We modify the natural $N$-equivariant 
structure on $M$ by multiplying by $\rho\prpr(n)$ 
the automorphism of $M$ given by $n$, for each $n\in N$. 
With this new structure, we have the identity $a(n)\circ a=a\circ n$ 
of automorphisms of (the total space of) $M$. To see this 
consider the identities
\begin{eqnarray*}\lefteqn{
(a^{-1}\circ(a(n).\rho\prpr(a(n)))\circ a)\circ(-n\rho\prpr(-n))}
\\&=&
(a^{-1}\circ a(n)\circ a\circ(-n)).\rho\prpr(a(n)).\rho\prpr(-n)
\\&=&
(a^{-1}\circ a(n)\circ a\circ(-n)).\rho(n)^{-1}={\rm Id}
\end{eqnarray*}
Consider now the bundle $(f_{*}M)^{N}$, which 
is the quotient of the bundle $M$ by the action 
of $N$. 
Using the identity $a(n)\circ a=a\circ n$, we see that the action 
of $a$ descends to 
$(f_{*}M)^N$. Furthermore, it is shown 
in \cite[Prop. 2, p. 70]{MAb} that $(f_{*}M)^N$ is 
naturally isomorphic to the total space of a line
bundle on $Y=X/N$. By 
 \cite[Prop. 2, p. 70]{MAb} again, this bundle has 
 the required properties.
\endProof
\begin{lemma}
Let $A$ be a $\Zn$-equivariant abelian variety over 
$\mtr{K}$, where $\Zn$ acts by isogenies. 
There exists a $\Zn$-equivariant abelian variety $B$ over $\mtr{K}$ and 
an ample $\Zn$-equivariant bundle $L$ on $B$, such that
\begin{description}
\item[(a)] $B$ is (non-equivariantly) isomorphic to 
$(A\times A^{\vee})^{4}$;  
\item[(b)] the equation $\chi(L)=1$ holds for 
the Euler characteristic of $L$;
\item[(c)] the group $\Zn$ acts by isogenies on $B$ and 
there exists an equivariant isogeny $p$ from 
$A^8$ to $B$.
\end{description}
\nlabel{eqzar}
\end{lemma}
\beginProof
We shall follow the steps of the proof of Zarhin's trick. 
Let $P\pr$ be any ample line bundle on $A$. 
Let $P:=\otimes_{g\in\Zn}g^{*}P\pr$. The bundle $P$ is ample
 and carries a natural $\Zn$-equivariant structure. 
Endow $A^4$ with the induced $\Zn$-equivariant structure and let 
$p_{i}:A^{4}\ra A$ be the $i$-th projection ($i=1,\dots 4$). Let
 $M\pr=p_{1}^{*}P\times p_{2}^{*}P\times p_{3}^{*} P\times p_{4}^{*}P$
be the fourth external tensor power of $P$. This 
is again an ample $\Zn$-equivariant line bundle on 
$A^4$. Let $m$ be the 
 order of the Mumford group $K(M\pr)$ of 
$M\pr$ (see \cite[p. 60]{MAb} for the definition). 
Now let $a,b,c,d$ be integers such that $a^2+b^2+c^2+d^2=-1(\mod\,m^2)$. Consider 
the endomorphism $\alpha(m)$ of $A^4$ described by the matrix
\begin{displaymath}
\left(\begin{array}{cccc}
a & -b & -c & -d\\
b & a & d & -c\\
c & -d & a & b\\
d & c & -b &a
\end{array}\right)
\end{displaymath}
(this is the endomorphism appearing in the proof of Zarhin's trick).  
This endomorphism commutes with all the elements of $\Zn$, because 
$\Zn$ acts by isogenies. 
Let now $N$ be the subgroup of $X=A^8$ given by the graph of 
$\alpha(m)|_{K(M\pr)}$. This subgroup is sent 
into itself by all the elements of $\Zn$. Let $B=X/N$ and let 
$p:X\ra B$ be the quotient map. By construction 
$B$ carries a natural $\Zn$-equivariant structure such that $p$ is 
equivariant. 
It is shown in \cite[Rem. 6.12, p. 136]{Milne} that 
there exists a line bundle
 $L$ on 
$B$ and an (non-equivariant) isomorphism 
 $M\pr\boxtimes M\pr\simeq p^*L$; we can 
 thus endow $M\pr\boxtimes M\pr$ with an $N$-equivariant 
 structure.  The $\Zn$-equivariant structure 
 on $M\pr\otimes_{\rm Ext}M\pr|_{X_{\Zn}}$ is by 
 construction trivial. By the last lemma, we 
can thus assume 
that $L$ carries a $\Zn$-equivariant structure and that 
there is an $\Zn$-equivariant isomorphism $M\pr\boxtimes M\pr\simeq p^*L$.
 We claim 
that $B$ and $L$ are the objects required in the 
statement of the lemma. The fact that (a) holds is 
a step in the proof of Zarhin's trick and 
we refer to \cite[Rem. 6.12, p. 136]{Milne} for the details. To see that 
(b) holds, 
we use \cite[Th. 2, p. 121]{MAb} and compute $\chi(M)=\chi(M\pr)^2=\#N$ and
$\chi(M)=\#N.\chi(L)$,   from which 
we deduce that $\chi(L)=1$. To see that (c) holds, note that 
the map $p$ defined above in the proof has the properties 
required of $p$ in (c). 
\endProof
We now quote 
the following results on extensions of line bundles from the generic fiber. If 
$X\ra S$ is any $S$ scheme, $L$ is a line sheaf over $X$ and $i:S\ra X$ is 
an $S$-valued point, then a {\it rigidification} of $L$ along $i$ is an 
isomorphism $i^*L\simeq {\cal O}_{S}$. The bundle $L$ together with 
a rigidification will be said to be {\it rigidified} 
along $i$. Once a point is given, 
the line bundles rigidified along that point form a category, 
where the morphisms are the sheaf morphisms that commute with the 
rigidification.
\begin{prop}
Let $A\ra\Spec K$ be an abelian variety over a number field. Suppose 
that $A$ has good reduction at all the finite places of $K$. Let 
${\cal A}\ra\Spec{\cal O}_{K}$ be its N\'eron model. 
Let $R_{\cal A}$ (resp. $R_{A}$) be the category of line 
bundles rigidified along the $0$-section 
of $\cal A$ (resp. $A$).
\begin{description}
\item[(1)](see \cite[1.1, p. 40]{MB}) The restriction functor $R_{\cal A}\ra R_{A}$ is an equivalence 
of categories;
\item[(2)](see \cite[Th. VIII.2]{R}) if $L$ is an ample line bundle on $A$, then any extension of $L$ to $\cal A$ is ample.
\end{description}
\nlabel{MBProp}
\end{prop}
Now let $K$ be a number field that contains all the $n$-th roots of unity. 
Let $A$ be an abelian variety over $K$ that has good reduction at 
all the finite places of $K$. Let $\cal A$ be its N\'eron 
model over ${\cal O}_K$ and let 
${\cal A}\pr:={\cal A}\times_{\Spec{\cal O}_{K}}
\Spec{\cal O}_{K}[{1\over n}]$. Let $f:{\cal A}\pr\ra\Spec{\cal O}_{K}[1/n]$ be 
the structure map and let ${\Omega}:=\Omega_{{\cal A}\pr}$
 be the sheaf of differentials of $f$.\\
Suppose that $A\pr$ is endowed 
with an action of $\Zn$ by $\Spec{\cal O}_{K}[{1\over n}]$-group 
scheme automorphisms. Let 
as usual $f^{\Zn}:{\cal A}\pr_{\Zn}\ra \Spec{\cal O}_{K}[{1\over n}]$ 
be the structure map of the fixed scheme. Let 
$i_{0}:\Spec{\cal O}_{K}[{1\over n}]\ra {\cal A}\pr$ be the zero 
section. We let $1\in\Zn$ act by $a(\cdot)$, as before.
\begin{lemma}
If $i^{*}_{0}{\Omega}$ has no ${\cal O}_{K}[1/n]$-submodule, which 
is fixed under the action of $\Zn$, then 
the scheme ${\cal A}\pr_{\Zn}$
is \'etale and finite over $\Spec{\cal O}_{K}[{1\over n}]$. 
\nlabel{finet}
\end{lemma}
\beginProof
We first prove that $f^{\Zn}$ is \'etale. We have to show that $f^{\Zn}$ is 
flat and that  
for any geometric point $\mtr{x}\ra\Spec{\cal O}_{K}[1/n]$, the corresponding 
scheme obtained by base-change is regular. The second condition is 
verified because of \cite[Cor. 2.9]{KR2} (base-change invariance of
 the fixed scheme) and \cite[Prop. 2.10]{KR2} 
 (regularity of the fixed-scheme). 
To show that $f^{\Zn}$ is flat, we apply the criterion 
\cite[Th. 9.9, p. 261]{Har}.
Choose a very ample $\Zn$-equivariant line bundle
$\cal L$ on 
$\cal A\pr$. Let $\frak p$ be a maximal ideal of $\Spec{\cal O}_{K}[1/n]$ and 
denote by $k({\frak p})$ the corresponding residue field.  
The Hilbert-Samuel polynomial of the fiber ${\cal A}\pr_{\Zn,k({\frak p})}$ 
of ${\cal A}\pr_\Zn$ at $\frak p$ relatively to 
 ${\cal
L}|_{{\cal A}\pr_{\Zn,k({\frak p})}}$  can be computed on the algebraic closure 
$\mtr{k({\frak p})}$ of
$k({\frak p})$. As 
$i_{0}^{*}{\Omega}$ has no fixed part, we see that 
${\cal A}\pr_{\Zn}(\mtr{k({\frak p})})$
consists  of isolated fixed points only. 
The number $F({\frak p})$ of these fixed
points is 
the Hilbert-Samuel polynomial of ${\cal A}_{\Zn,k({\frak p})}\pr$, which 
is of degree $0$. We have to show that $F({\frak p})$ is independent 
of $\frak p$. To prove this, consider first 
that $H^{1}({\cal A}\pr_{k({\frak p})},{\cal L}_{k({\frak p})})=0$ 
for all $\frak p$; 
this follows from the characterization of the cohomology 
of line bundles on abelian varieties \cite[Par. 16, p. 150]{MAb}. Thus 
 we know that $f_{*}{\cal L}$ is locally free and that 
 there are natural equivariant isomorphisms 
 $(f_{*}{\cal L})_{k({\frak p})}
 \simeq f_{k({\frak p}),*}{\cal L}_{k({\frak p})}$.
  We may also assume 
that the action on ${\cal L}|_{{\cal A}\pr_{\Zn}}$ is trivial; this might 
be achieved by replacing $\cal L$ by its $n$-th 
tensor power. Using the Lefschetz trace formula of 
\cite{BFQ}, we see that 
$F({\frak p})$ depends only on the trace of 
$a$ on $H^{0}({\cal A}_{K},{\cal L}_{K})$ and on the 
determinant of ${\rm Id}-a$ on $i_{0,K}^{*}\Omega_{K}$ 
(where $(\cdot)_{K}$ refers to the base change by $\Spec\ K\ra 
\Spec{\cal O}_{K}[{1\over n}])$. 
Thus 
$F({\frak p})$ is independent of $\frak p$. 
To see that $f^{\Zn}$ is finite, we only have to check that it 
is quasi-finite, as it is projective (see \cite[Ex. 11.2]{Har}). Let again 
$\frak p$ be a prime ideal. As $f^{\Zn}$ is 
\'etale, we know that ${\cal A}\pr_{\Zn,k({\frak p})}$ is the spectrum 
of a direct sum of finite field extensions of $k({\frak p})$; furthermore 
this sum is finite, since the morphism is of finite type. 
Hence ${\cal A}\pr_{\Zn,k({\frak p})}$ is a finite set and thus we are done.
\endProof 

\section{The equivariant analytic torsion of line bundles 
on abelian varieties}

Let $(V,g^V)$ be a d-dimensional Hermitian vector space and let $\Lambda\subset V$ be
a lattice of rank $2 d$. The quotient $A:=V/\Lambda$ is a flat complex torus.
According to the Appell-Humbert theorem
\cite[Ch. 2.2]{LB}, the holomorphic line bundles on $A$ can be described as follows:
Choose an Hermitian form $H$ on $V$ such that $E:={\rm Im}\,H$ takes integer values
on $\Lambda\times \Lambda$. Choose furthermore $\alpha:\Lambda\to S^1$ such that
$$
\alpha(\lambda_1+\lambda_2)=\alpha(\lambda_1)\alpha(\lambda_2) e^{\pi i
E(\lambda_1,\lambda_2)}
$$
for all $\lambda_1,\lambda_2\in\Lambda$. Then there is an associated line bundle
$L_{H,\alpha}$ defined as the quotient of the trivial line bundle on $V$ by the action
of $\Lambda$ given by
$$
\lambda \circ (v,t):=(v+\lambda,\alpha(\lambda) e^{\pi H(v,\lambda)+\pi
H(\lambda,\lambda)/2}t)\,\,.
$$
There is a canonical Hermitian metric $h^L$ on $L_{H,\alpha}$ given by
$$
h^L((v,t_1),(v,t_2)):= t_1\bar t_2 e^{-\pi H(v,v)}\,\,.
$$
Define $C\in{\rm End}(V)$ by $H(v_1,v_2)=g^V(v_1, C v_2)$. $C$ is Hermitian with
respect to $g^V$. Consider an automorphism $g$ of $(V,g^V)$, leaving $\Lambda$, $H$
and
$\alpha$ invariant. Then $g$ and $C$ commute, thus they may be diagonalized
simultaneously. Denote their eigenvalues by $(e^{i \phi_j})$ and $(\nu_j)$,
respectively, and let
$(e_j)$ be a corresponding set of $g^V$-orthonormal eigenvectors. Assume that for all
$j$,
$\phi_j\notin 2\pi{\bf Z}$, i.e. that $g$ acts on $A$ with isolated fixed points. Note
that the isometric automorphism
$g$ of
$A$ has finite order.

Let $L_{\rm Tr}(L_{H,\alpha})$ denote the trace of the action of $g$ on $H^0(A,L_{H,\alpha})$.
\begin{lemma}\nlabel{zeromu}
Assume that $\nu_j=0$ for all $j$. Then the
equivariant analytic torsion of
$\bar L_{H,\alpha}$ on
$(A, g^V)$ with respect to $g$ vanishes.
\end{lemma}
\beginProof
Let $V^\vee_{\bf R}$ be the dual of the underlying real vector space of $V$. Consider
the dual lattice $\Lambda^\vee:=\{\mu\in V^\vee_{\bf R}|\mu(\lambda)\in 2\pi {\bf
Z}\forall\lambda\in\Lambda\}$. Represent $\alpha=e^{i \mu_0}$ with $\mu_0\in
V^\vee_{\bf R}$. It is shown in 
\cite[section V]{Koehler4} 
that the
eigenfunctions are given by the functions $f_\mu:V_{\bf R}/\Lambda\to{\bf C}$,
$x\mapsto e^{i(\mu+\mu_0)(x)}$ with corresponding eigenvalue
$\frac1{2}\|\mu+\mu_0\|^2$. The eigenforms are given by the product of these
eigenfunctions times the pullback of elements of $\Lambda^*V^\vee$. The eigenvalue of
such a form $f_\mu\cdot \eta$ is the eigenvalue of $f_\mu$.

As $L_{H,\alpha}$ is $g$-invariant, we get $g^*e^{i \mu_0}=e^{i (\mu_0+\mu_1)}$ for
some $\mu_1\in\Lambda^\vee$. Also,
$g$ maps $\Lambda^\vee$ to itself, thus for any $\mu\in\Lambda^\vee$ there is some
$\mu'\in\Lambda^\vee$ with $g^*f_\mu=f_{\mu'}$.
 As $g$ acts fixed point free on
$V\setminus\{0\}$, it maps a function
$f_\mu$ to a multiple of itself iff $\mu+\mu_0=0$, i.e. iff $f_\mu$ represents an
element in the cohomology. Thus $g$ acts diagonal free on the complement of the
cohomology, and the zeta function defining the torsion vanishes.
\endProof

Let $\psi:{\bf R}/{\bf Z}\to {\bf R}$ denote the function $\psi([x]):=(\log
\Gamma)'(x)$ for $x\in]0,1]$.

\begin{theor} Assume that $\nu_j>0$ for all $j$.
Then the
equivariant analytic torsion of
$\bar L_{H,\alpha}$ on
$(A, g^V)$ with respect to $g$ is given by
\begin{eqnarray*}
T_g(A,\bar L_{H,\alpha})&=&\sum_{j=1}^d\Bigg[\frac{\log(2\pi \nu_j)}{e^{-i
\phi_j}-1} +i R^{\rm rot}(\phi_j)
\\&&
-\frac1{4}
\left(2\log(2\pi)-2\Gamma'(1)+\psi([\frac{\phi_j}{2\pi}])+\psi([1-\frac{\phi_j}{2\pi}])\right)
\Bigg]L_{\rm Tr}(L_{H,\alpha})\,\,.
\end{eqnarray*}
\nlabel{KaiTh}
\end{theor}
\beginProof
Assume first that the $\nu_j$ are pairwise linear independent over $\bf Q$. 
Then as is shown in \cite[p. 3]{Ber}, the spectrum of the Kodaira-Laplace operator on
$A$ is given by
$$
\sigma(\square)=\big\{2\pi\sum_{j=1}^d n_j \nu_j\,\big|n_j\in{\bf N}_0\ \forall
j\big\}\,\,.
$$
Set $e^i:=g^V(e_i,\cdot)$. Let $E^q_\nu$ denote the eigenspace corresponding to the
eigenvalue $\nu$ of $\square$ on $\Gamma^\infty(A,\Lambda^q T^*A\otimes
L_{H,\alpha})$. We shall prove by induction that the trace of $g$ on $E^q_{2\pi\sum_j
n_j|\nu_j|}$ is given by
\begin{equation}\nlabel{eigform}
{\#\{j|n_j\neq0\}\choose q}L_{\rm Tr}(L_{H,\alpha})\cdot\prod_j e^{i n_j \phi_j}\,\,.
\end{equation}
First this formula holds for $E^0_0$.
Now the eigenspaces verify the relations
\begin{equation}\nlabel{ber1}
E^q_\nu=\bigoplus_{j_1<\dots<j_q}\left(\Lambda_{k=1}^q e^{j_k}\otimes
E^0_{\nu-2\pi\sum_{k=1}^q \nu_{j_k}}\right)
\end{equation}
\cite[eq. (7)]{Ber} and the complex
\begin{equation}\nlabel{ber2}
0\to E^0_\nu\stackrel{\bar\partial}\to\cdots \stackrel{\bar\partial}\to E^d_\nu\to0
\end{equation}
is acyclic. Thus we can determine the trace on $E^q_\nu$ for $q>0$ by the trace on
some $E^0_{\nu'}$ with $\nu'<\nu$. Then the trace on $E^0_\nu$ is given by the trace
on the $E^q_\nu$ for $q>1$ by the sequence \ref{ber2}. As the relations
(\ref{ber1}), (\ref{ber2}) are compatible with (\ref{eigform}), equation
(\ref{eigform}) is proven.

Hence the zeta function defining the torsion is given by
$$
Z(s)=\sum_{q=0}^d \sum_{\nu\in\sigma(\square)} \frac{(-1)^{q+1}q}{\nu^s}
{\#\{j|n_j\neq0\}\choose q}L_{\rm Tr}(L_{H,\alpha})\cdot\prod_j e^{i n_j \phi_j}\,\,.
$$
Notice that for any $0\leq k\leq n$
$$
\sum_{q=0}^d (-1)^{q+1}q{k\choose q}=\left\{
\begin{array}{ccc}
1&\mbox{ for } k=1&\\
0&\mbox{ for } k\neq1&\,\,.
\end{array}\right.
$$
Consider for $z\in{\bf C}$, $|z|=1$ and for $s\in{\bf C}$, Re $s>1$ the zeta function
$$
L(z,s):=\sum_{k=1}^\infty \frac{z^k}{k^s}
$$
and its meromorphic continuation to $s\in{\bf C}$.
Thus
$$
Z(s)=\sum_{j=1}^d \sum_{k=1}^\infty \frac{e^{i k \phi_j}}{(2\pi k\nu_j)^s} L_{\rm
Tr}(L_{H,\alpha}) =\sum_{j=1}^d (2\pi \nu_j)^{-s}L(e^{i \phi_j},s) L_{\rm
Tr}(L_{H,\alpha})
$$
and
$$
Z'(0)=\sum_{j=1}^d \left(-\log(2\pi \nu_j)L(e^{i \phi_j},0)+L'(e^{i \phi_j},0)
\right)L_{\rm Tr}(L_{H,\alpha})\,\,.
$$
Using the definition of $R^{\rm rot}(\phi):=\frac{\partial}{\partial s}_{|s=0}
\frac1{2 i}(L(e^{i \phi},s)-L(e^{-i \phi},s))$ in
\cite{Koehler} (compare section 2), the formula for
$\frac{\partial}{\partial s}_{|s=0}
\frac1{2}(L(e^{i \phi},s)+L(e^{-i \phi},s))$ in
\cite[Lemma 13]{Koehler3} and the formula $L(e^{i \phi},0)=(e^{-i\phi}-1)^{-1}$ (e.g. in
\cite[p. 108]{Koehler2}) we find the theorem for $\nu_j$ pairwise linear independent
over $\bf Q$. If the
cohomology does not change, the torsion varies continuously with the metric
(\cite[section (d)]{BGS3}),
thus the result holds for any nonzero $\nu_j$.
Alternatively, one can show that formulas (\ref{ber1}), (\ref{ber2}) hold more
general by arguing as in
\cite[Remark p. 4]{Ber}) or by continuity.
\endProof
{\bf Remark.} More general, assume only that all $\nu_j$ are non-zero. Again, by the
results of \cite{Ber} and a proof similar to the above one we get
\begin{eqnarray*}
T_g(A,\bar L_{H,\alpha})&=&\sum_{j=1}^d {\rm sign}(\nu_j) \Bigg[\frac{\log(2\pi
|\nu_j|)}{e^{-i
\phi_j}-1} +i R^{\rm rot}(\phi_j)
\\&&
-\frac1{4}
\left(2\log(2\pi)-2\Gamma'(1)+\psi([\frac{\phi_j}{2\pi}])+\psi([1-\frac{\phi_j}{2\pi}])\right)
\Bigg]L_{\rm Tr}(L_{H,\alpha})\,\,.
\end{eqnarray*}
Combining this with Lemma \ref{zeromu} as in \cite[section 4]{Ber} by splitting $A$
(see also \cite[chapter 3,\S3]{LB}) and using the product formula for equivariant
torsion \cite[Lemma 2]{Koehler2}, one obtains the value of the equivariant torsion for
any
$L_{H,\alpha}$. Similarly, for automorphisms $g$ having a larger fixed point set one
can split $A$ accordingly to obtain the value of the torsion.

\section{Application of the fixed point formula to abelian varieties
with C.M. by a cyclotomic field}

Let now $n>0$ be a natural number and 
let $\phi(n):=
\#({\bf Z}/(n))^{\times}$.  
Let $A$ be an abelian variety 
of dimension $d=\phi(n)/2$ defined 
over a number field $K$. Suppose that $A_{\mtr{\bf Q}}$ has complex
multiplication by ${\cal O}_{{\bf Q}(\mn)}$ and 
fix a ring embedding ${\cal O}_{{\bf Q}(\mn)}\ra 
\End(A_{\mtr{\bf Q}})$ 
(see \cite{Sh} for the general theory). As before 
we choose a primitive $n$-th root of unity 
$\zeta:=\zeta_{n}$; this root defines an isomorphism 
$\mn({\bf C})\simeq\Zn$ and thus a canonical 
$\Zn$-action on $A_{\mtr{\bf Q}}$. 
We may suppose that 
${\bf Q}(\mn)\subseteq K$ and that $A$ has good reduction 
at all the finite places of $K$. The latter hypothesis 
is possible in view of \cite[Th. 7, p. 505]{ST}. 
Let $B$ be the abelian variety obtained from $A_{\mtr{\bf Q}}$ 
via lemma \ref{eqzar}. We can suppose without loss of generality that $B$, 
as well as the line bundle 
$L$ promised in \ref{eqzar}, are defined over $K$; we may also 
suppose that the map $p$ appearing in (c) is defined over $K$.
The existence of $p$ shows that $B$ has 
good reduction at the finite places of $K$ as well.  
Furthermore, we normalize the 
action of $\Zn$ on $L$ so that its restriction to $L|_0$ becomes 
trivial (one might achieve this by multiplying 
the action by some character of $\Zn$). Choose 
some rigidification of $L$ along the $0$-section. The action of $\Zn$ will 
then leave this rigidification 
invariant.\\
We let $\pi:{\cal
B}\ra\Spec{\cal O}_K$ be  the N\'eron model of $B$. 
By the universal property of N\'eron models, the action 
of $\Zn$ on $B$ extends to $\cal B$. Also, by (1) of Proposition 
\ref{MBProp},  
there exists a line bundle $\cal L$ on $\cal B$, endowed with 
a $\Zn$-action and a $\Zn$-invariant rigidification along 
the $0$-section, which both extend the corresponding structures on $L$. 
By \cite[Prop. 2.1, p. 48]{MB}, there exists a unique 
metric on ${\cal L}_{\Bbb C}$ whose first Chern form 
is translation invariant and such that the rigidification is an isometry. 
We endow $L_{\bf C}$ with this metric. It is $\Zn$-invariant by unicity. 
We endow $B({\bf C})$ with the translation invariant K\"ahler metric 
given by the Riemann form of $L_{\Bbb C}$. We let 
${\Omega_\pi}$ be the sheaf of differentials of $\pi$ and endow 
it with the metric induced by the K\"ahler metric. If we 
let $i_{0}:\Spec{\cal O}_{K}\ra{\cal B}$ be the $0$-section, we then 
have 
a canonical isometric equivariant isomorphism 
$\pi^*i_{0}^*\mtr{{\Omega_\pi}}\simeq\mtr{{\Omega_\pi}}$. 
It is shown in \cite[Par. 4, 4.1.10]{BAb} that 
$\ar{\deg}(i_{0}^*\mtr{{\Omega_\pi}})=[K:{\Bbb
Q}]h_{\rm Fal}(B_{\Qb})$.\\   
We now apply the equivariant arithmetic Riemann-Roch formula 
\refeq{RBCEq} to $\mtr{\cal L}$. We shall work over 
${\cal O}_{K}[1/n]$, over which $\Zn$- and $\mn$-action 
are equivalent concepts because of lemma \ref{CNMN}. Let $d_{B}=8d$ be 
the relative  dimension of $\cal B$. Let $p_{1},\dots ,p_{l}$ be 
the prime factors of $n$ and let $U$ be the complement of 
the set $\{p_{1},\dots,p_{l}\}$ in $\Spec{\Bbb Z}$. We then 
have the identification $U_{K}\simeq \Spec{\cal O}_{K}[1/n]$. 
We identify $\mn$ and $\Zn$ over $U_{K}$ via the root of unity $\zeta$.  
We let ${\cal S}$ be the subgroup of $\Bbb R$ generated by
 the expressions
$\sum_{k=1}^{l}q_{k}{\rm Im}(\alpha_{k})\log p_{k}$, where 
$q_{k}\in{\Bbb Q}$ is a rational number and 
$\alpha_{k}$ is an $n$-th 
root of unity.
\begin{prop}
Let $e^{i\phi_{1}},\dots e^{i\phi_{d_{B}}}$ be the 
eigenvalues of the action of $1\in\Zn$ on 
$TB|_{0}$ (with multiplicities). The equality
\begin{equation}
{1\over [K:{\bf Q}]}\sum_{k\in\Zn}
\NIm{{\zeta^{k}\over 1-\zeta^k}}\ar{\deg}(i_0^*\mtr{\Omega_{\pi,k}})=
\frac1{2}\sum_{j=1}^{d_{B}}R^{\rm rot}(\phi_j)
\nlabel{ModH}
\end{equation}
\nlabel{EARRAb}
holds in ${\bf R}/{\cal S}$. 
\end{prop}
\beginProof
For any $l\in{\Bbb Z}$, let $[l]$ denote the $l$-plication map 
from a group scheme into itself. 
First notice that by \cite[(4.1.23)]{BAb} (i.e. a hermitian 
refinement of the theorem of the square), the identity
$[l]^*\ac1(\mtr{\cal L}\otimes[-1]^*\mtr{\cal L})=
l^2\ac1(\mtr{\cal L}\otimes[-1]^*\mtr{\cal L})$ holds 
in $\ar{\rm CH}^{1}({\cal B})$
, for 
every $l\in{\Bbb Z}$. Now consider that by lemma \ref{finet}, 
the scheme ${\cal B}_\mn$ is
 a finite commutative group scheme over ${\cal O}_{K}[1/n]$; 
 it thus has an 
order $l_0\geq 0$ and the $l_{0}$-plication map is the $0$-map on 
${\cal B}_\mn$ (see \cite{T}). Now we compute in $\ar{\rm CH}^{1}({\cal B}_\mn)$:
$$
[l_{0}]^*\ac1(\mtr{\cal L}|_{{\cal B}_\mn}\otimes[-1]^*\mtr{\cal L}|_{{\cal B}_\mn})
=
l_{0}^2\ac1(\mtr{\cal L}|_{{\cal B}_\mn}\otimes[-1]^*\mtr{\cal L}|_{{\cal B}_\mn})=0
$$
and thus we get 
$\ar{\deg}\blb f\umn_{*}(l_{0}^2\ac1(\mtr{\cal L}|_{{\cal B}_\mn})+l_{0}^2
\ac1([-1]^*\mtr{\cal L}|_{{\cal B}_\mn}))\brb=
\ar{\deg}\blb f\umn_{*}(2.l_{0}^2\ac1(\mtr{\cal L}|_{{\cal B}_\mn}))\brb=0$ and 
finally $\ar{\deg}\blb f\umn_{*}(\ac1(\mtr{\cal L}_{{\cal B}
_\mn}))\brb=0$.\\
Let us now write down the equation \refeq{RBCEq} in our situation. 
First notice that in view of the equivariant isomorphism 
$\pi^*i_{0}^*{{\Omega_\pi}}\simeq{{\Omega_\pi}}$, the functions 
$T$ and $D$ are constant. 
Until the end of the proof, let $N$ be the number of fixed points 
of the action of $1$ on an arbitrary connected component 
of $B({\bf C})$. The holomorphic Lefschetz trace 
formula (see \cite{BFQ} or \cite[III, (4.6), p. 566]{ASe}) shows 
that $N.T/D=\sum_{k\in\Zn}\zeta_{n}^{k}\rk((R^{0}f_{*}{\cal L})_{k})$ 
and thus $N$ is independent of the choice of the component. 
We obtain
\begin{eqnarray*}\lefteqn{
\sum_{k\in{\Bbb Z}/(n)}\zeta^k .\ar{\deg}((R^{0}f_{*}\mtr{\cal L})_k)=
}
\\&=&
{1\over 2}T_{g}(B({\Bbb C}),\mtr{E})-
i[K:{\bf Q}]{N.T\over D}.\sum_{k=1}^{d_{B}}
\rrot(\phi_{k})
\\&+&
\ar{\deg}\blb f\umn_{*}({T\over D}.
\sum_{k\in\Zn}{\zeta^{k}\over 1-\zeta^{k}}
\ac1(\mtr{\Omega}_{\pi,k}))\brb
\end{eqnarray*}
in ${\bf C}/{\cal C}$. As $\rk(R^{0}\pi_{*}({\cal L}))=1$, we also have 
the equality of complex numbers
$$
\sum_{k\in{\Bbb Z}/(n)}{\zeta^k}.\ar{\deg}((R^{0}f_{*}\mtr{\cal L})_k)=
\ar{\deg}((R^{0}f_{*}\mtr{\cal L})).
{N.T\over D}
$$
And thus we get
\begin{eqnarray}\lefteqn{
\ar{\deg}((R^{0}f_{*}\mtr{\cal L}))
}\nonumber\\&=&
({N.T\over D})^{-1}{1\over 2}T_{g}(B({\Bbb C}),\mtr{E})-
i[K:{\bf Q}]\sum_{k=1}^{d_{B}}
\rrot(\phi_{k})\nonumber\\
&+&
\ar{\deg}\blb f\umn_{*}(\sum_{k\in\Zn}{\zeta^{k}\over (1-\zeta^{k})}
\ac1(\mtr{\Omega}_{\pi,k}))\brb.N^{-1}\label{YAEq}
\end{eqnarray}
in ${\bf C}/{\cal C}$. 
The theorem \ref{KaiTh} shows that the imaginary part 
of $({N.T\over D})^{-1}{1\over 2}T_{g}(B({\Bbb C}),\mtr{E})$ is 
equal to 
$[K:{\bf Q}]\rrot(\phi_{k})$. 
Furthermore, using the 
isomorphism $\pi^*i_{0}^*\mtr{{\Omega_\pi}}\simeq\mtr{{\Omega_\pi}}$ and 
the projection formula, we see that
$$ 
\ar{\deg}\blb f\umn_{*}(\sum_{k\in\Zn}{\zeta^{k}\over 1-\zeta^{k}}
\ac1(\mtr{\Omega}_{\pi,k}))\brb.N^{-1}=
\sum_{k\in\Zn}{\zeta^{k}\over 1-\zeta^{k}}
\ac1(i_{0}^{*}\mtr{\Omega}_{\pi,k}).
$$
Taking these two facts into account, we can take the imaginary 
part of both sides in \refeq{YAEq} to conclude.
\endProof
Let now $\Phi:=
\{\Phi_{1},\dots,\Phi_{d}\}$ be the type of $A_{\mtr{\bf Q}}$. 
Define $\zeta^{P(k)}=\Phi_{k}(\zeta)$ for $k=1,\dots ,d$. 
In view of \ref{eqzar} (c), the identity \refeq{ModH} can be rewritten as
$$
{1\over[K:{\bf Q}]}\sum_{k=1}^{d}
\NIm{{\Phi_{k}(\zeta)\over 1-\Phi_{k}(\zeta)}}
\ar{\rm deg}(i_{0}^{*}\mtr{\Omega}_{\pi,P(k)})=
-4\sum_{k=1}^{d} {\partial\over\partial s}{L}^{\rm Im}(\Phi_{k}(\zeta),0)
$$
(in ${\bf R}/{\cal S}$). 
Now notice that we can change our choice of
$\zeta$ in the latter equation and replace it by 
$\sigma(\zeta)$, where $\sigma\in{\rm Gal}(\Qmn|{\bf Q})$
(this corresponds to applying the fixed point formula to 
a power of the original automorphism), thus 
obtaining a system of linear equations in the 
$\ar{\rm deg}(i_{0}^{*}\mtr{\Omega}_{\pi,k})$. Notice 
also the equations obtained from $\sigma(\zeta)$ and
$\mtr{\sigma(\zeta)}$ are equivalent. With this remark 
in mind, we see that the just mentioned system of equations 
is equivalent to the following one:
\begin{equation}
{1\over[K:{\bf Q}]}\sum_{k=1}^{d}
\NIm{{\Phi_{k}^{-1}\circ\Phi_{l}(\zeta)\over 1-\Phi_{k}^{-1}
\circ\Phi_{l}(\zeta)}}X_{k}=
-4\sum_{k=1}^{d} {\partial\over\partial s}{L}^{\rm Im}
(\Phi_{k}^{-1}\circ\Phi_{l}(\zeta),0)+E_{l}
\nlabel{MainEq}
\end{equation}
(in $\bf R$) where $l=1,\dots,d$ and the coefficients of the vector 
$E:={^t}[E_{1},\dots,E_{d}]$ lie in $\cal S$. We shall show 
that the matrix $M:=[\NIm{{\Phi_{k}^{-1}
\circ\Phi_{l}(\zeta^{l})\over 1-\Phi_{k}^{-1}
\circ\Phi_{l}(\zeta)}}]_{l,k}$ of this system is invertible 
as a matrix of real numbers. Since the coefficients of $M$ all lie 
in $\Qmn$, the coefficients of $M^{-1}$ lie in $\Qmn$ as well and 
thus we see that the coefficients of the vector $M^{-1}E$ lie 
in $\cal C$. Thus we can determine the quantities 
$X_{k}$ up to an element of $\cal S$. 
Now recall that by construction
$$
h_{\rm Fal}({\cal B}_{\Qb})={1\over [K:{\bf Q}]}(
\sum_{k=1}^{d}\ar{\rm deg}(i_{0}^{*}\mtr{\Omega}_{\pi,P(k)})).
$$
in ${\bf R}/{\cal S}$. Furthermore, by
a result of Raynaud in 
 \cite[Exp. VII, Cor. 2.1.3, p. 207]{SCM}, the modular 
 height of an abelian variety and the modular height 
 of the dual of the latter are equal. Thus, 
 using (a) \ref{eqzar}, we see that
$h_{\rm Fal}({\cal B}_{\mtr{\Bbb Q}})=8.h_{\rm Fal}(A_{\mtr{\bf Q}})$ (in $\bf R$). 
Thus we can determine $h_{\rm Fal}(A_{\mtr{\bf Q}})$ up to an element
of $\cal S$.\\
Before we proceed to solve the system \refeq{MainEq}, 
we need two lemmata that relate the quantities appearing 
in the system with Dirichlet $L$-functions.\\
{\bf Warning}. In what follows, in contradiction with classical 
usage, the notation $L(\chi,s)$ will 
always refer to the {\bf non-primitive} $L$-function associated
with a Dirichlet character. We write $\chi_{\rm prim}$ for 
the primitive character associated with $\chi$ and accordingly 
write $L(\chi_{\rm prim},s)$ for the associated primitive 
$L$-function.\\
First take notice of the elementary fact that 
${\rm Im}(z/(1-z))={1\over 2}\cot({1\over 2}\arg(z))$ if 
$|z|=1$. Let $G$ be the Galois group of the extension 
${\Bbb Q}(\mu_{m})|{\Bbb Q}$. 
From now on, for simplicity we fix $\zeta=e^{2\pi i/n}$. 
The following lemma is a variation on the functional equation 
of the Dirichlet $L$-functions. When $\chi$ is a primitive 
character, it can be derived directly from the functional 
equation and classical results on Gauss sums. 
\begin{lemma}
Let $\chi$ be an odd character of $G$. 
The equality
$$
<{L}^{\rm Im}(\sigma(\zeta),s),\chi>:=
{1\over 2d}\sum_{\sigma\in G}\mtr{\chi(\sigma)}
{L}^{\rm Im}(\sigma(\zeta),s)
={1\over 2d}
n^{1-s}{\Gamma(1-s/2)\over\Gamma((s+1)/2)}\pi^{s-1/2}
L(\mtr{\chi},1-s)
$$
holds for all $s\in {\bf C}$. 
\nlabel{SixMonthLem}
\end{lemma}
In the expression $<{L}^{\rm Im}(\sigma(\zeta),s),\chi>$, 
the symbol $\sigma$ is considered as a variable in $G$. 
\beginProof
We prove the equality for $0<s<1$. The full equality 
then follows by analytic continuation. We compute
\begin{eqnarray}\lefteqn{
n^{{s+1\over 2}}\Gamma({s+1\over 2})\pi^{-{s+1\over 2}}{1\over 2d}
\sum_{\sigma\in G}\mtr{\chi(\sigma)}{L}^{\rm Im}(\sigma(\zeta),s)}\nonumber
\\&=&
n^{{s+1\over 2}}\Gamma({s+1\over 2})\pi^{-{s+1\over 2}}{1\over 2d}
\sum_{\sigma\in G}\sum_{k\geq 0}{{\rm Im}(\sigma(\zeta)^{k})\over k^{s}}
\mtr{\chi}(\sigma)\nonumber
\\&=&
n^{{s+1\over 2}}\Gamma({s+1\over 2})\pi^{-{s+1\over 2}}{1\over 2d}
\sum_{k\geq 0}{1\over k^{s}}\sum_{\sigma\in G}\mtr{\chi}(\sigma)
{\rm Im}(\sigma(\zeta)^{k})\nonumber
\\&=&
\sum_{k\geq 0}\int_{0}^{\infty}{1\over 2d}
k.e^{-k^{2}\pi u/n}u^{{1\over 2}(s+1)-1}[\sum_{\sigma\in G}
\mtr{\chi}(\sigma){\rm Im}(\sigma(\zeta)^{k})]{\rm\bf du}\nonumber
\\&=&
-i.\int_{0}^{\infty}{1\over 2d}\sum_{k\geq 0}
k.e^{-k^{2}\pi\over nu}u^{-1-{1\over 2}(s+1)}[\sum_{\sigma\in G}
\mtr{\chi}(\sigma)\sigma(\zeta)^{k}]{\rm\bf du}\nlabel{LebEq}
\\&=&
\int_{0}^{\infty}{1\over 2d}\sum_{k\geq 0}
\chi(k)u^{-1-{1\over 2}(s+1)}{k\over n}
e^{-\pi k^{2}u/n}(n.u)^{3\over 2}{\rm\bf du}\nlabel{PoissonEq}
\\&=&
{1\over 2d}\sum_{k\geq 0}
\chi(k)\int_{0}^{\infty}u^{-1-{1\over 2}(s+1)}{k\over n}
e^{-\pi k^{2}u/n}(n.u)^{3\over 2}{\rm\bf du}\nlabel{LCruxEq}
\\&=&
{1\over 2d}
n^{3/2-s/2}\Gamma(1-s/2)\pi^{s/2-1}L(\mtr{\chi},1-s)\nonumber
\end{eqnarray}
For the equality \refeq{PoissonEq}, we used the Poisson 
summation formula. To obtain the equality \refeq{LebEq}, 
we first exchange the summation and integration symbols and then 
make the change of variable $u\mapsto{1\over u}$. The exchange 
of summation and integration symbols is justified
by the following estimates. Let 
$\sigma_0\in G$ and let $t_{k}:={\rm Im}(\sigma_{0}(\zeta)^{k})$ and 
$v_{k}:=k.e^{k^{2}\pi u/n}$. Since $\sum_{k=0}^{ln}t_{k}=0$ 
for all $l\geq 0$, the sequence $T_{k}:=
\sum_{j=0}^{k}t_{j}$ is bounded above by $C>0$. Consider now 
that 
\begin{equation}
|\sum_{k=0}^{N}k.e^{-k^{2}\pi u/n}{\rm Im}(\sigma_{0}(\zeta)^{k})|=
|\sum_{k=0}^{N}t_{k}v_{k}|=
|\sum_{k=0}^{N-1}T_{k}(v_{k}-v_{k+1})+T_{N}v_{N}|\nlabel{PartSumEq}
\end{equation}
where we used partial summation for the last equality. 
The function (of $k$) $k.e^{-k^{2}\pi u/n}$ is increasing  on the 
interval $[0,\sqrt{n\over  2\pi u}]$ and decreasing 
on the interval $[\sqrt{n\over 2\pi u},\infty[$. Let 
$k_0$ be the largest integer less or equal to 
$\sqrt{n\over  2\pi u}$. 
The expression \refeq{PartSumEq} can be bounded above by 
\begin{eqnarray*}\lefteqn{
C v_{N}+C\sum_{k=k_{0}+1}^{N-1}(v_{k}-v_{k+1})+
C|v_{k_{0}}-v_{k_{0}+1}|+
C\sum_{k=0}^{k_{0}-1}(v_{k+1}-v_{k})}
\\&=&
Cv_{N}+C(v_{k_{0}+1}-v_{N})+C(v_{k_0}-v_0)+
C|v_{k_{0}}-v_{k_{0}+1}|=
C[v_{k_{0}+1}+v_{k_{0}}+|v_{k_{0}}-v_{k_{0}+1}|]
\\&=&
2C\max\{v_{k_{0}+1},v_{k_{0}}\}\leq 
2C\sqrt{n\over 2\pi e}u^{-1/2}.
\end{eqnarray*}
Hence 
$$
|\sum_{k=0}^{N}k.e^{-k^{2}\pi u/n}u^{{1\over 2}(s+1)-1}
{\rm Im}(\sigma_{0}(\zeta)^{k})|\leq 
2C\sqrt{n\over 2\pi e}|u^{-1/2}u^{{1\over 2}(s+1)-1}|=
2C\sqrt{n\over 2\pi e}|u^{{1\over 2}s-1}|
$$
On the other hand, for $u>1$, the classical estimate
$$
|\sum_{k=0}^{N}k.e^{-k^{2}\pi u/n}u^{{1\over 2}(s+1)-1}
{\rm Im}(\sigma_{0}(\zeta)^{k})|\leq 
|\sum_{k=0}^{N}k.e^{-k\pi u/n}u^{{1\over 2}(s+1)-1})|\leq 
$$
$$
\leq e^{-\pi u/n}|u^{{1\over 2}(s+1)-1}|{1\over (1-e^{-\pi u/n})^{2}}
$$
holds. The first estimate show that 
for $u\in]0,1[$ and $s>0$, the function of $u$
$$|\sum_{k=0}^{N}k.e^{-k^{2}\pi u/n}u^{{1\over 2}(s+1)-1}
{\rm Im}(\sigma_{0}(\zeta)^{k})|$$ is bounded by 
an element of $L^{1}([0,1])$. The second estimate 
shows that for $u\in]1,\infty[$ and  
$s\in {\bf C}$ the same function is bounded by 
 an element of $L^{1}(]1,\infty[)$. 
By the dominated convergence theorem, we might thus exchange 
summation and integration 
symbols. Completly similar estimates justify the equality \refeq{LCruxEq}.  
\endProof
It is shown in \cite[Proof of Th. 8, p. 564]{Koehler} that 
${L}^{\rm Im}(\sigma(\zeta),0)=
{1\over 2}\cot({1\over 2}\arg(\sigma(\zeta)))$. Thus, we see that 
the last lemma has the following corollary. 
\begin{cor}
For any odd character $\chi$ on $G$, the equality 
$$
\sum_{\sigma\in G}\cot({1\over 2}\arg(\sigma(\zeta)))\chi(\sigma)=
{2n\over \pi}L(\chi,1)
$$
holds. 
\nlabel{ColLem}
\end{cor} 
To tackle with the system \refeq{MainEq},
 we shall need the following result from 
linear algebra. We say that a complex-valued 
function $f$ on $G$ is {\it odd} if 
$f(c\circ x)=-f(x)$ for all $x\in G$, where $c$ denotes
complex conjugation. 
\begin{lemma}
Let $\vec{X}:=(X_{1},\dots, X_{d})$ and 
$\vec{Y}:=(Y_{1},\dots ,Y_{d})$. Let $f$ be any odd 
function on $G$ and let $M_f$ be the matrix 
$[f(\Phi_{k}^{-1}\circ\Phi_{l})]_{l,k}$. If $<\chi,f>\not = 0 $ 
for all the odd characters $\chi$, then the system $M\vec{X}=\vec{Y}$ 
is maximal and 
$$
X_{j}=\sum_{\chi\ {\rm odd}}{\chi(\Phi_{j})\sum_{l}\mtr{\chi}(\Phi_{l})Y_{l}
\over
d^{2}<f,\chi>}
$$
\nlabel{WasLem}
\end{lemma}
\beginProof
Let $\phi_{k}$ be the function defined on $G$ such that 
$\phi_{k}(x)=1$ if $x=\Phi_{1}\circ\Phi_{k}$ and $0$ otherwise ($k=1,\dots,d$). 
on $G$. Let $V^-$ be the complex vector space of 
odd functions on $G$. An ordered basis $B_\Phi$ of 
$V^-$ is given by $\phi_1-\phi_1\circ c,\dots, 
\phi_d-\phi_d\circ c$. Another basis $B_\chi$ is 
given by the odd characters on $G$. We view $M$ as a linear 
endomorphism of $V^-$, via $B_\Phi$. We now proceed to find the matrix 
of $M$ in the basis $B_\chi$. We compute 
$$
M.^t[\chi(\Phi_1\circ\Phi_1),\dots,\chi(\Phi_1\circ\Phi_d)]=
^t[\sum_k f(\Phi_k^{-1}\circ\Phi_l)\chi(\Phi_1\circ\Phi_k)]_l
$$
and 
\begin{eqnarray*}
\lefteqn{
\sum_k f(\Phi_k^{-1}\circ\Phi_l)\chi(\Phi_1\circ\Phi_k)=
{1\over 2}\sum_{\sigma\in G} f(\sigma^{-1}\circ\Phi_l)\chi(\Phi_1\circ\sigma)}
\\&=&
{1\over 2}\sum_{\sigma\in G} f(\sigma^{-1})
\chi(\sigma\circ\Phi_{l}\circ\Phi_{1})=
d<f,\chi>\chi(\Phi_{1}\circ\Phi_{l}).
\end{eqnarray*}
Thus $M$ is represented by the diagonal matrix
${\rm Diag}[d<f,\chi>]_\chi$ in the basis $B_\chi$. 
The vector $\vec{Y}_\chi:=
{1\over d}[\sum_l Y_k\mtr{\chi}(\Phi_1\circ\Phi_l)]_\chi$ 
represents the vector $\vec{Y}$ in $B_\chi$. Thus 
the solution of $M\vec{X}=\vec{Y}$ in $B_\chi$ is the vector 
$\vec{X}_\chi:=[{\sum_l Y_k\mtr{\chi}(\Phi_1\circ\Phi_l)\over 
d^{2}<f,\chi>}]_\chi$ and we thus obtain 
\begin{equation}
X_{j}=\sum_{\chi\ {\rm odd}}\chi(\Phi_{1}\circ\Phi_{j}){
\sum_l Y_k\mtr{\chi}(\Phi_1\circ\Phi_l)\over 
d^{2}<f,\chi>}=
\sum_{\chi\ {\rm odd}}{\chi(\Phi_{j})\sum_{l}\mtr{\chi}(\Phi_{l})Y_{l}
\over
d^{2}<f,\chi>}
\end{equation}
\endProof
In the following proposition, we apply the lemma \ref{WasLem}
 to the system 
\refeq{MainEq} and use the lemmata \ref{ColLem} and 
\ref{SixMonthLem} to evaluate the resulting 
expression in terms of logarithmic derivatives of $L$-functions. 
\begin{prop}
The Faltings height of $A_{\mtr{\bf Q}}$ is given by 
the identity
$$
{1\over d}h_{\rm Fal}(A)=-\sum_{\chi\ {\rm odd}}<\Phi*\Phi^{\vee},\chi>
2{L\pr(\chi_{\rm prim},0)\over L(\chi_{\rm prim},0)}+\sum_{p|n}a_{p}\log(p)
$$
where $a_{p}\in\Qmn$.
\nlabel{MainProp}
\end{prop}
\beginProof
We apply the lemma \ref{WasLem} to the system \refeq{MainEq} 
(with $E_{l}=0$). 
Define $f:G\ra {\bf C}$ by the 
formula $f(\sigma):=\cot({1\over 2}\arg(\sigma(\zeta)))$. 
The fact that the system is maximal is implied by the 
fact that $L(\chi,1)\not =0$, when $\chi$ is a non-principal 
Dirichlet character (see for instance 
\cite[Th. 2, p. 212]{CaFr}). We compute
\begin{eqnarray*}\lefteqn{
{1\over 8}\sum_j X_j = -\sum_j\sum_{\chi\ {\rm odd}} {\chi(\Phi_j)\sum_{k,l}
{\partial\over \partial s}{L}^{\rm Im}
(\Phi_k^{-1}\circ\Phi_l(\zeta),0)\mtr{\chi}(\Phi_l)
\over d^{2}<f,\chi>}}
\\&=&
-\sum_\chi{\sum_j\chi(\Phi_j)\sum_{k,l}
{\partial\over \partial s}{L}^{\rm Im}(\Phi_k^{-1}
\circ\Phi_l(\zeta),0)\mtr{\chi}(\Phi_l)\over
d^{2}<f,\chi>}
\end{eqnarray*}
Using scalar and convolution products, we can write
\begin{eqnarray}\lefteqn{
{\sum_j\chi(\Phi_j)\sum_{k,l}
{\partial\over \partial s}{L}^{\rm Im}(\Phi_k^{-1}\circ\Phi_l(\zeta),0)
\mtr{\chi}(\Phi_l)\over
d^{2}<f,\chi>}}\nonumber
\\&=&
{2d<\chi,\Phi>.2d^{2}<\chi,{\partial\over \partial s}
{L}^{\rm Im}(\sigma(\zeta),0)*\Phi*(\cdot)^{-1}>
\over 
d^{2}<f,\chi>}\nonumber
\\&=&
{4d^{3}<\chi,\Phi><\chi,{\partial\over \partial s}
{L}^{\rm Im}(\sigma^{-1}(\zeta),0)*\Phi^{\vee}>
\over d^{2}<f,\chi>}\nonumber
\\&=&
8d^{2}<\chi,\Phi*\Phi^{\vee}>.{
<{\partial\over \partial s}{L}^{\rm Im}
(\sigma(\zeta),0),\chi>\over {2n\over\pi}L(\mtr{\chi},1)}\nonumber
\\&=&
{4d^{2}\pi\over n}<\Phi*\Phi^{\vee},\chi>.
{<{\partial\over \partial s}{L}^{\rm Im}(\sigma(\zeta),0),\chi>\over 
L(\mtr{\chi},1)}\nlabel{FLEq}
\end{eqnarray}
where $\sigma\in G$ is a variable and the equalities 
are in $\bf R$.  
Now by lemma \ref{SixMonthLem}, the sum over all odd $\chi$ of \refeq{FLEq} 
is equal to
\begin{equation}
\sum_{\chi\ {\rm odd}}<\Phi*\Phi^{\vee},\chi>(\ {4d^2 \pi\over n}{{\partial\over\partial s}[
L(\mtr{\chi},1-s)n^{1-s}{\Gamma(1-s/2)\over \Gamma((s+1)/2)}\pi^{s-1/2}]_{s=0}{1\over 
2d}\over L(\mtr{\chi},1)}\ )
\nlabel{IDKWEq}
\end{equation}
Furthermore, from the existence of 
Euler product expansions for $L$ functions, 
we deduce that 
$$
{L\pr(\mtr{\chi},1)\over 
L(\mtr{\chi},1)}={L\pr(\mtr{\chi}_{\rm prim},1)\over 
L(\mtr{\chi}_{\rm prim},1)}+\sum_{p|n}{\mtr{\chi}_{\rm prim}(p)/p\over 
1-\mtr{\chi}_{\rm prim}(p)/p}\log(p).
$$ 
The sum
$$
\sum_{\chi\ {\rm odd}}<\Phi*\Phi^{\vee},\chi>
{\mtr{\chi}_{\rm prim}(p)/p\over 
1-\mtr{\chi}_{\rm prim}(p)/p}
$$
is an algebraic number. By construction, it is invariant 
under the action of any element of ${\rm Gal}(\mtr{\bf Q}|{\bf Q})$ and 
it is thus an element of $\bf Q$. Taking this into account 
and using the functional equation of primitive Dirichlet 
$L$-functions, we obtain that \refeq{IDKWEq} is equal to 
\begin{eqnarray}\lefteqn{
\sum_{\chi\ {\rm odd}}<\Phi*\Phi^{\vee},\chi>(\ {2d\pi\over n}[{-n\Gamma\pr(1)\over 2\pi}-{n\log(n)\over\pi}+
{n\log(\pi)\over \pi}-{n\Gamma\pr(1/2)\over 2\Gamma(1/2)\pi}]}\nonumber
\\&+&
2d[{L\pr(\chi_{\rm prim},0)\over L(\chi_{\rm prim},0)}-
\log({2\pi\over f_{\mtr{\chi}}})+{\Gamma\pr(1)\over \Gamma(1)}]\ )
\nlabel{ALEq}
\end{eqnarray}
in ${\bf R}/{\cal S}$. By Galois invariance again, the expressions 
$\sum_{\chi\ {\rm odd}}<\Phi*\Phi^{\vee},\chi>\log(n)$ and 
$\sum_{\chi\ {\rm odd}}<\Phi*\Phi^{\vee},\chi>\log(f_{\mtr{\chi}})$ lie 
in $\cal S$. Furthermore, by \cite[(1), p. 363]{C2}, 
we have the equality $\sum_{\chi\ {\rm odd}}<\Phi*\Phi^{\vee},\chi>={1\over 2}$ 
of real numbers. Thus the two lines of the expression \refeq{ALEq} are 
equal to 
\begin{eqnarray}\lefteqn{
[2d\sum_{\chi\ {\rm odd}}<\Phi*\Phi^{\vee},\chi>
{L\pr(\chi_{\rm prim},0)\over L(\chi_{\rm prim},0)}]+
2d[-\Gamma\pr(1)/4+{1\over 2}\log(\pi)-{1\over 4}{\Gamma\pr(1/2)\over \Gamma(1/2)}
-{1\over 2}\log(2\pi)+{1\over 2}{\Gamma\pr(1)\over\Gamma(1)}]}\nonumber
\\&=&
[2d\sum_{\chi\ {\rm odd}}<\Phi*\Phi^{\vee},\chi>
{L\pr(\chi_{\rm prim},0)\over L(\chi_{\rm prim},0)}]
\\&+&
2d[-\Gamma\pr(1)/4-{1\over 2}\log(2)+{1\over 2}{\Gamma\pr(1)\over\Gamma(1)}-
{1\over 4}{\Gamma\pr(1)\over\Gamma(1)}+{1\over 2}\log(2)]\nlabel{GammaEq}
\\&=&
2d\sum_{\chi\ {\rm odd}}<\Phi*\Phi^{\vee},\chi>
{L\pr(\chi_{\rm prim},0)\over L(\chi_{\rm prim},0)}\nonumber
\end{eqnarray}
in ${\bf R}/{\cal S}$.  
We used the equality $2^{2z-1}\Gamma(z)\Gamma(z+{1\over 2})=\sqrt\pi
\Gamma(2z)$ ($z\in{\bf C}$) for the equality \refeq{GammaEq}. 
This concludes the proof. 
\endProof
The next lemma shows that the expression for the Faltings 
height appearing in the last proposition is invariant under 
number field extension (this lemma is implicit in the
definition of $h_{\rm Fal}$ given in \cite{C2}). 
\begin{lemma}
Let $E\subseteq E\pr$ be C.M. fields. Suppose that 
$E$ and $E\pr$ are both Galois extensions of $\bf Q$. 
Let $\Phi$ be a type of $E$ and let $\Phi\pr$ the type induced 
by $\Phi$ on $E\pr$. The identity of real 
numbers
$$
\sum_{\chi_{E}}<\Phi*\Phi^{\vee},\chi_{E}>
{L\pr(\chi_{E,{\rm prim}},0)\over L(\chi_{E,{\rm prim}},0)}=
\sum_{\chi_{E\pr}}<\Phi\pr*{\Phi\pr}^{\vee},\chi_{E\pr}>
{L\pr(\chi_{E\pr,{\rm prim}},0)\over L(\chi_{E\pr,{\rm prim}},0)}
$$
holds, where the first sum is over the odd characters $\chi_{E}$ of 
$E$ and the second sum is over the odd characters of 
$E\pr$.
\end{lemma}
\beginProof
Let ${\cal G}_{E}$ (resp. ${\cal G}_{E\pr}$) be the 
Galois group of $E$ (resp. of $E\pr$) over $\bf Q$. 
Let $p:{\cal G}_{E\pr}\ra {\cal G}_{E}$ be the natural map. 
By construction the identities
$$ 
<f\circ p,g\circ p>_{E\pr}=<f,g>_{E}\ {\rm and} 
\ (f\circ p)*(g\circ p)=(f*g)\circ p
$$ 
hold for functions 
$f,g:{\cal G}_{E}\ra {\bf C}$. By construction again, 
the equality $L((\chi\circ p)_{\rm prim},s)=
L(\chi_{\rm prim},s)$ also holds, where $\chi$ is a character 
of ${\cal G}_{E}$. To prove the 
identity of the lemma, it is thus sufficient to show that 
$<\Phi\pr*{\Phi\pr}^{\vee},\chi>=0$ if $\chi$ is an  
odd character of $E\pr$ not induced from $E$.  
So let $\chi$ be such a character. By assumption, 
$\chi$ is non-trivial on the normal 
subgroup $H:=p^{-1}({\rm Id}_{E})$ of ${\cal G}_{E\pr}$ and 
$(\Phi\pr*{\Phi\pr}^{\vee})(h.g)=
(\Phi\pr*{\Phi\pr}^{\vee})(g)$ for all $g\in{\cal G}_{E\pr}$ and 
$h\in H$. Let $h_{0}\in H$ be such that $\chi(h_{0})\not= 1$. 
We compute 
\begin{eqnarray*}\lefteqn{
<\Phi\pr*{\Phi\pr}^{\vee},\chi>:=
{1\over\#{\cal G}_{E\pr}}\sum_{g\in{\cal G}_{E\pr}}
(\Phi\pr*{\Phi\pr}^{\vee})(g)\mtr{\chi}(g)}
\\&=&
{1\over\#{\cal G}_{E\pr}}\sum_{g\in{\cal G}_{E\pr}}
(\Phi\pr*{\Phi\pr}^{\vee})(g.h_{0})\mtr{\chi}(g.h_{0})
=[{1\over\#{\cal G}_{E\pr}}\sum_{g\in{\cal G}_{E\pr}}
(\Phi\pr*{\Phi\pr}^{\vee})(g)\mtr{\chi}(g)]\mtr{\chi}(h_{0})
\\&=&
<\Phi\pr*{\Phi\pr}^{\vee},\chi>\mtr{\chi}(h_{0})
\end{eqnarray*}
If $<\Phi\pr*{\Phi\pr}^{\vee},\chi>\not= 0$, then 
$\mtr{\chi}(h_{0})=1$, a contradiction.
\endProof
A slightly weaker form of the following lemma (in the sense that 
only the isogeny classes of the varieties are determined) is 
contained in a lemma due to Shimura-Taniyama. 
\begin{lemma}
Let $E\subseteq\mtr{\bf Q}$ be a C.M. field that is an 
abelian extension of $\bf Q$ 
and let $\Phi$ be a type 
of $E$. Fix an integer $n$ such 
that $E\subseteq{\bf Q}(\mn)$ and let $\Phi\pr$ be the type 
of ${\bf Q}(\mn)$ lifted from $\Phi$ via the inclusion. Let 
$A$ be an abelian variety over $\mtr{\bf Q}$, that admits 
a complex multiplication by ${\cal O}_{E}$ of type $\Phi$.
 The abelian variety 
$A^{[{\bf Q}(\mn):E]}$ admits a complex multiplication 
by ${\cal O}_{{\bf Q}(\mn)}$ of type $\Phi\pr$. 
\end{lemma}
\beginProof
Let $r:=[{\bf Q}(\mn):E]$. 
From the fact that ${\cal O}_{{\bf Q}(\mn)}$ is generated as 
an $E$-algebra by a primitive 
root of $1$, we deduce that there exists a basis of $r$ elements 
of ${\cal O}_{{\bf Q}(\mn)}$ as an ${\cal O}_{E}$-module. Notice 
now that multiplication by a fixed element of ${\cal O}_{{\bf Q}(\mn)}$ 
defines an ${\cal O}_{E}$-module endomorphism of 
${\cal O}_{{\bf Q}(\mn)}$, that is the identity iff 
the element is $1$. 
We can thus use the 
basis to define a ring injection ${\cal O}_{{\bf Q}(\mn)}
\ra M_{r\times r}({\cal O}_{E})$ of ${\cal O}_{{\bf Q}(\mn)}$ 
into the ring of $r\times r$ matrices with coefficients in 
${\cal O}_{E}$. This ring injection gives rise to 
a complex multiplication by ${\cal O}_{{\bf Q}(\mn)}$ on the product 
$A^{r}$. By construction, each element 
$\sigma\in{\rm Hom}({\bf Q}(\mn),\mtr{\bf Q})$ 
of the corresponding type has the property that 
$\sigma|_{E}\in\Phi$. This implies that the type is lifted from 
$\Phi$ and concludes the proof.
\endProof
\begin{cor}
Let $A$ be an abelian variety of dimension $d$ defined 
over $\mtr{\bf Q}$. Let $E$ be a C.M. field, suppose 
that $E$ is an abelian extension of $\bf Q$ and that 
$A$ has complex multiplication by ${\cal O}_{E}$. Let
${\cal O}_{E}\ra End(A)$ be an embedding of rings and let  
$\Phi$ be the associated type. Then 
the identity
$$
{1\over d}h_{\rm Fal}(A)=-\sum_{\chi\ {\rm odd}}<\Phi*\Phi^{\vee},\chi>
2{L\pr(\chi_{\rm prim},0)\over L(\chi_{\rm prim},0)}+\sum_{p|f}a_{p}\log(p)
$$
of real numbers holds, where $f$ is the conductor of $E$ over $\bf Q$ and 
$a_{p}\in{\bf Q}(\mu_{f})$.
\nlabel{MainCor}
\end{cor}
\beginProof
By class field theory $E\subseteq{\bf Q}(\mu_f)$; furthermore 
the modular height $h_{\rm Fal}$ is additive for products of 
abelian varieties over $\mtr{\bf Q}$. We can thus 
apply the two last lemmata to reduce the proof of the 
identity to the case $E={\bf Q}(\mu_f)$. This case is covered 
by proposition \ref{MainProp} and so we are done.
\endProof
{\bf Remark}. The equation \refeq{MainEq} can also be deduced from the 
application of the arithmetic fixed point formula to the 
trivial bundle, instead of an ample bundle. To perform 
this alternative computation, one has to 
use the computation of the relative cohomology of the trivial 
bundle on an abelian scheme given in \cite{BBM}. This 
cohomology can then be related to the sheaf of 
relative differentials via the equivariant principal polarisation 
provided by section 2.

\end{document}